\documentclass[10pt,a4paper]{article}
\usepackage[textheight=236mm,textwidth=151mm]{geometry}
\usepackage{amsthm,amssymb,graphicx}

\title{\bf Space-Efficient Karatsuba Multiplication for Multi-Precision Integers}
\author{Yiping Cheng\\ School of Electronic and Information Engineering\\ Beijing Jiaotong University, Beijing 100044\\{\tt ypcheng@bjtu.edu.cn}}
\date{}

\begin{document}

\maketitle

\begin{abstract}
The traditional Karatsuba algorithm for the multiplication of
polynomials and multi-precision integers has a time complexity of
$O(n^{1.59})$ and a space complexity of $O(n)$. Roche proposed an
improved algorithm with the same $O(n^{1.59})$ time complexity but
with a much reduced $O(\log n)$ space complexity. In Roche's paper
details were provided for multiplication of polynomials, but not
for multi-precision integers. Multi-precision integers differ from
polynomials by the presence of carries, which poses difficulties
in implementing Roche's scheme in multi-precision integers. This
paper provides a detailed solution to these difficulties. Finally,
numerical comparisons between the schoolbook, traditional
Karatsuba, and space-efficient Karatsuba algorithms are provided.
\end{abstract}

\section{Introduction}

In primary schools pupils are taught how to multiply two integers.
This algorithm has been implemented in the multiplication of
multi-precision integers (or called long integers), known as the
schoolbook algorithm. It requires $O(n^2)$ time and $O(1)$ space.
Here $n$ is the length, i.e., the number of limbs in the
representation of the two equal-length integers for
multiplication. Karatsuba and Ofman in \cite{karatsuba63} noticed
that using a divide-and-conquer technique, an $O(n^{1.59})$
sub-quadratic time complexity can be achieved, with the price
being an additional $O(n)$ memory space. Aiming to remove the need
of extra space, Roche proposed in \cite{roche09} a space-efficient
variant of the Karatsuba algorithm. However, the description of
the algorithm was mainly addressed to the multiplication of
polynomials. The implementation of space-efficient Karatsuba in
multi-precision integers is not straightforward, as there are
carries when adding or subtracting integers, even when we use a
subtractive version of Karatsuba. Tricky issues resulting from the
presence of carries need to be resolved.

This paper will give a detailed implementation of space-efficient
Karatsuba for multiplication of multi-precision integers. We
follow the basic scheme of \cite{roche09}, but a few changes have
been made, as listed below: \begin{itemize} \item A subtractive
version is adopted instead of the additive version.

\item Four separate variables for storing the carries are
introduced, because the results of additions and subtractions may
fall beyond the original representation.

\item For even length square multiplication, three cases are
identified and are processed for obtaining correct results in all
possible cases. This subdivision is not needed in polynomial
multiplication but is unavoidable in long integer multiplication.

\item A different treatment for odd length is used, which is
chosen for more efficient carry processing.
\end{itemize}

The rest of the paper is organized as follows. Section 2 gives the
mathematical foundation for the space-efficient Karatsuba for long
integers. Section 3 describes our C implementation of this
algorithm. Section 4 provides numerical results for performance
comparison. We conclude the paper with Section 5.

\section{Mathematical Formulation}
\subsection{Representation of Long Integers}

Long integers are represented by {\em limbs\/}. Let $\rho$ be the
radix. Then any integer $A\in \{0..\rho^{n}-1\}$ can be
represented uniquely with $n$ ordered limbs, each limb $a_i$ being
an integer within $\{0..\rho-1\}$. That is,

\begin{equation} A=\sum_{i=0}^{n-1} a_i\rho^{n-1-i}.\end{equation}

For intermediate values in the algorithm, it is sometimes useful
to store it with a carry. If $A$ falls (not very far) beyond the
range $\{0..\rho^{n}-1\}$, then we can represent it with $n$ limbs
and a carry $\kappa$, which can be positive, zero, or negative.
That is

\begin{equation} A=\kappa \rho^n + \sum_{i=0}^{n-1} a_i\rho^{n-1-i}.\end{equation}

In computers we can use an {\tt int\/} variable to store the
carry, and as the range of a 32-bit {\tt int\/} variable is
$\{-2^{31}..2^{31}-1\}$, it will be more than sufficient for our
purposes here.

\subsection{Even Equal Length Additive Multiplication}

As in \cite{roche09}, we consider the following problem, which we
call {\em equal length additive multiplication problem\/} here.
Given four multi-precision integers $A^{(0)},A^{(1)},B,C\in
\{0..\rho^{n}-1\}$, to compute
\begin{equation} D=(A^{(0)}-A^{(1)})B+C\rho^{n}.\end{equation} In
this subsection we deal with the case when $n$ is even, i.e.,
$n=2k$.

Let
\begin{eqnarray}
A^{(0)} &=& A^{(0)}_0\rho^{k}+A^{(0)}_1, \\
A^{(1)} &=& A^{(1)}_0\rho^{k}+A^{(1)}_1,\\
B &=& B_0\rho^{k}+B_1, \\
C &=& C_0\rho^{k}+C_1,
\end{eqnarray}
where $A^{(0)}_0,A^{(0)}_1,A^{(1)}_0,A^{(1)}_1,B_0,B_1,C_0,C_1\in
\{0..\rho^{k}-1\}$.

We then have \[ (A^{(0)}-A^{(1)})B= \]\begin{equation}
(A^{(0)}_0-A^{(1)}_0)B_0\rho^{2k}+[(A^{(0)}_0-A^{(1)}_0)B_1+(A^{(0)}_1-A^{(1)}_1)B_0]\rho^{k}
+ (A^{(0)}_1-A^{(1)}_1)B_1.
\end{equation}

Now define the following intermediate multi-precision integers:
\begin{eqnarray}
P_{0} &=& (A^{(0)}_0-A^{(1)}_0)B_0, \\
P_{1} &=& (A^{(0)}_1-A^{(1)}_1)B_1, \\
E &=& A^{(0)}_0-A^{(1)}_0-A^{(0)}_1+A^{(1)}_1,\\
P_2 & = & (B_1-B_0)E.
\end{eqnarray}

It is easy to see that \[ -2(\rho^k-1)\leq E \leq 2(\rho^k-1).\]

For $E$ we divide three cases:

\subsubsection{Case 1: $-(\rho^k-1)\leq E\leq \rho^k-1$}

This case can be decided when the internal $k$ limbs and one carry
representation of $E$ has a carry $0$ or $-1$, but $E\neq
-\rho^k$. How do we know whether $E=-\rho^k$? We note that
$E=-\rho^k$ if and only if it has carry $-1$ and all $0$ limbs.

In this case, if $E\geq 0$, we can use the original equation (12)
to compute $P_2$, otherwise, use the mathematically equivalent
equation
\begin{equation} P_2 = (B_0-B_1)(-E).
\end{equation}

By (8), we have \begin{equation}
(A^{(0)}-A^{(1)})B=P_{0}\rho^{2k}+(P_{0}+P_{1}+P_2)\rho^{k} +
P_{1}.
\end{equation}

Now suppose \begin{eqnarray}
P_{0} &=& P_{00}\rho^{k}+P_{01}, \\
P_{1} &=& P_{10}\rho^{k}+P_{11}, \\
P_2 &=& P_{20}\rho^{k}+P_{21}.
\end{eqnarray} where $P_{01},P_{11},P_{21}\in
\{0..\rho^{k}-1\}$.

Then \begin{equation}
(A^{(0)}-A^{(1)})B=P_{00}\rho^{3k}+(P_{01}+P_{10}+P_{00}+P_{20})\rho^{2k}+(P_{01}+P_{10}+P_{11}+P_{21})\rho^{k}
+P_{11}.
\end{equation}

The above four subsequences of limbs can be computed using a
procedure very similar to \cite[Section 2.2]{roche09}, except that
\begin{itemize}
\item As here we use a subtractive version, we need to be careful
about the signs.

\item Four {\tt int\/} variables
$\kappa_{00},\kappa_{01},\kappa_{10},\kappa_{10}$ are introduced,
to store the carries resulting from addition and subtraction.
\end{itemize}

The computation steps should be clear from Table 1. Note that
recursive calls take place in steps 3,5,8.

\begin{table}[h]
\centering

\caption{Snapshot of $D$ and the $\kappa$'s through the steps for
case 1}

\begin{tabular}{c|cccc}
 & $\kappa_{00}\rho^k+D_{00}$ & $\kappa_{01}\rho^k+D_{01}$ & $\kappa_{10}\rho^k+D_{10}$ & $\kappa_{11}\rho^k+D_{11}$ \\
\hline
0 & $C_0$ & $C_1$ \\
1 & $C_0$ & $C_1-C_0$ \\
2 & $C_0$ & $C_1-C_0$ & &  $|E|$\\
3 & $C_0$ & $C_1-C_0+P_{20}$ & $P_{21}$  \\
4 & $C_0$ & $C_1-C_0+P_{20}$ & $P_{21}$ & $C_1-C_0+P_{20}-P_{21}$ \\
5 & $C_0+P_{00}$ & $P_{01}$ & $P_{21}$ & $C_1-C_0+P_{20}-P_{21}$ \\
6 & $C_0+P_{00}$ & $P_{01}$ & $P_{01}+P_{21}$ & $C_1-C_0+P_{20}-P_{21}$ \\
7 & $C_0+P_{00}$ & $C_1+P_{00}-P_{21}+P_{20}$ & $P_{01}+P_{21}$ \\
8 & $C_0+P_{00}$ & $C_1+P_{00}-P_{21}+P_{20}$ & $P_{01}+P_{21}+P_{10}$ & $P_{11}$ \\
9 & $C_0+P_{00}$ & $C_1+P_{01}+P_{10}+P_{00}+P_{20}$
& $P_{01}+P_{21}+P_{10}$ & $P_{11}$ \\
10 & $C_0+P_{00}$ & $C_1+P_{01}+P_{10}+P_{00}+P_{20}$ &
$P_{01}+P_{21}+P_{10}+P_{11}$ & $P_{11}$
\end{tabular}
\end{table}

\subsubsection{Case 2: $\rho^k\leq E\leq 2(\rho^k-1)$}

In this case, $P_2$ falls beyond the range $\{0..\rho^k-1\}$, and
the number corresponding to the pure limbs representation of $E$
is $E-\rho^k$. So we compute
\begin{equation} Q_2 = (B_1-B_0)(E-\rho^k).\end{equation}

Thus \[ P_2 = Q_2+(B_1-B_0)\rho^k.\]

And we have
\begin{eqnarray}
P_{20} &=& Q_{20}+B_1-B_0, \\
P_{21} &=& Q_{21}.
\end{eqnarray} And \[(A^{(0)}-A^{(1)})B= \]
\begin{equation}
P_{00}\rho^{3k}+(P_{01}+P_{10}+P_{00}+\underbrace{Q_{20}+B_1-B_0}_{=P_{20}})\rho^{2k}+(P_{01}+P_{10}+P_{11}+\underbrace{Q_{21}}_{=P_{21}})\rho^{k}
+P_{11}.
\end{equation}

In this case the computation steps should be clear from Table 2.
\begin{table}[h]
\centering

\caption{Snapshot of $D$ and the $\kappa$'s through the steps for
case 2}

\begin{tabular}{c|cccc}
 & $\kappa_{00}\rho^k+D_{00}$ & $\kappa_{01}\rho^k+D_{01}$ & $\kappa_{10}\rho^k+D_{10}$ & $\kappa_{11}\rho^k+D_{11}$ \\
\hline
0 & $C_0$ & $C_1$ \\
1 & $C_0$ & $C_1-C_0$ \\
2 & $C_0$ & $C_1-C_0$ & &  $E-\rho^k$\\
3 & $C_0$ & $C_1-C_0+Q_{20}$ & $P_{21}$  \\
4 & $C_0$ & $C_1-C_0+\underbrace{Q_{20}+B_1-B_0}_{=P_{20}}$ & $P_{21}$  \\
5 & $C_0$ & $C_1-C_0+P_{20}$ & $P_{21}$ & $C_1-C_0+P_{20}-P_{21}$ \\
6 & $C_0+P_{00}$ & $P_{01}$ & $P_{21}$ & $C_1-C_0+P_{20}-P_{21}$ \\
7 & $C_0+P_{00}$ & $P_{01}$ & $P_{01}+P_{21}$ & $C_1-C_0+P_{20}-P_{21}$ \\
8 & $C_0+P_{00}$ & $C_1+P_{00}-P_{21}+P_{20}$ & $P_{01}+P_{21}$ \\
9 & $C_0+P_{00}$ & $C_1+P_{00}-P_{21}+P_{20}$ & $P_{01}+P_{21}+P_{10}$ & $P_{11}$ \\
10 & $C_0+P_{00}$ & $C_1+P_{01}+P_{10}+P_{00}+P_{20}$
& $P_{01}+P_{21}+P_{10}$ & $P_{11}$ \\
11 & $C_0+P_{00}$ & $C_1+P_{01}+P_{10}+P_{00}+P_{20}$ &
$P_{01}+P_{21}+P_{10}+P_{11}$ & $P_{11}$ \\
\end{tabular}
\end{table}

\subsubsection{Case 3: $-2(\rho^k-1)\leq E\leq -\rho^k $}

In this case, we compute \begin{equation}  Q_2 =
(B_0-B_1)(-E-\rho^k). \end{equation}

Then \[ P_2 = Q_2+(B_0-B_1)\rho^k.\]
\begin{eqnarray}
P_{20} &=& Q_{20}+B_0-B_1, \\
P_{21} &=& Q_{21}.
\end{eqnarray}

The corresponding equation and table are similar to case 2, with
$E$ replaced with $-E$ and $B_1-B_0$ replaced with $B_0-B_1$. So
they will not be reproduced here.

\subsubsection{Final Carry Processing}

It should be noted that when the steps in Tables 1 and 2 are
completed, $\kappa_{11}=0$ is automatically ensured, and then we
have \[ D=\kappa_{00}\rho^{4k} + D_{00}\rho^{3k}+
\kappa_{01}\rho^{3k}+
D_{01}\rho^{2k}+\kappa_{10}\rho^{2k}+D_{10}\rho^{k}+D_{11}.\]

We then need to do a final step of carry processing, so that
finally $\kappa_{01}=\kappa_{10}=0$, and

\[ D=\kappa_{00}\rho^{4k} + D_{00}\rho^{3k}+D_{01}\rho^{2k}+D_{10}\rho^{k}+D_{11}.\]

\subsection{Odd Equal Length Additive Multiplication}

Now consider the same problem as given in the previous subsection,
but now with $n$ odd, i.e. $n=2k+1$.

There are a number of possible expansions that can reduce this
case to the even case. But we found the following expansion most
convenient:

Let \begin{eqnarray}
A^{(0)}-A^{(1)} &=& \rho (\bar{A}^{(0)}-\bar{A}^{(1)}) + a^{(0)}- a^{(1)}\\
B &=& \rho^{2k} b + \bar{B} \\
C &=& \rho^{2k} c + \bar{C}
\end{eqnarray} where $\bar{A}^{(0)},\bar{A}^{(1)},\bar{B},\bar{C}\in \{0..\rho^{2k}-1\}$ and $a^{(0)},a^{(1)},b,c\in
\{0..\rho-1\}$.

We have
\[ (A^{(0)}-A^{(1)})B+C\rho^{2k+1}= \rho^{2k}(a^{(0)}- a^{(1)})b + (a^{(0)}- a^{(1)})\bar{B} + \]
\begin{equation}
\rho^{2k+1}(\bar{A}^{(0)}-\bar{A}^{(1)})b+ \rho
(\bar{A}^{(0)}-\bar{A}^{(1)})\bar{B}
+\rho^{4k+1}c+\rho^{2k+1}\bar{C}.\end{equation}

We then have
\[ (A^{(0)}-A^{(1)})B+C\rho^{2k+1}= \rho [(\bar{A}^{(0)}-\bar{A}^{(1)})\bar{B}
+\rho^{2k}\bar{C}] +\rho^{4k+1}c\]
\begin{equation}
+\rho^{2k}(A^{(0)}- A^{(1)})b + (a^{(0)}- a^{(1)})\bar{B}
.\end{equation}

On the right-hand side of (30), the first term has already been
computed using the algorithm of the previous section, and the
second term is already present in the buffer to store the product,
and only the latter two terms need to be computed and added to the
product now. The advantage of this expansion is that the two
remaining terms involve nonintersecting subsequences of the entire
limbs sequence of the product, and therefore the carries can thus
be processed most efficiently.

\subsection{General Length Multiplication}

We now come back to our original problem: Suppose $A$ is a long
integer of length $n$, and $B$ is a long integer of length $m$. We
are to compute
\[AB.\]

Let $n\geq m\geq 1$, and suppose
\begin{equation}
n = qm+r, \mbox{ where } 0\leq r<m.
\end{equation}

Then we let \begin{equation}A = A_0\rho^{qm} +
A_1\rho^{(q-1)m}+\cdots+A_{q-1}\rho^m + A_q.\end{equation} where
$A_0\in\{0..\rho^r-1\}$, and $A_1,\cdots,A_q\in\{0..\rho^m-1\}$.

This general length multiplication is done by a recursive routine.
If $r=0$, then no recursive call of itself is needed, but we need
to fill the first $m$ limbs of the product buffer to zero.
Otherwise, the routine first calls itself to compute $A_0B$, which
is stored in the leading $m+r$ limbs of the product buffer. Let

\[ A_0B = (A_0B)_0\rho^m + (A_0B)_1\]
 where $(A_0B)_0\in\{0..\rho^r-1\}$, and
$(A_0B)_1\in\{0..\rho^m-1\}$.

Then, we call the equal length additive routine to compute

\[ A_1B + (A_0B)_1\rho^m,\] and process the carry, so that \[A_0B\rho^m + A_1B\] is computed and stored in the buffer.

Go on with these steps repeatedly until the entire $AB$ is
computed. This requires $q$ calls of the equal length additive
routine, and no heap memory allocation is needed here.

\section{Implementation of Algorithm in C}

The algorithm mathematically described above has been implemented
in C. The implementation consists of two groups of functions: the
primitive group, and the upper level group.

\subsection{The Primitive Group}

The primitive group of functions include the following functions:

{\tt int MpiAdd(uint $\rho$, LIMB* A, const LIMB* B, int n);}

{\tt int MpiSub(uint $\rho$, LIMB* A, const LIMB* B, int n);}

{\tt int MpiNeg(uint $\rho$, LIMB* A, int n);}

{\tt int MpiAddC(uint $\rho$, LIMB* A, int n, int $\kappa$);}

Here the type {\tt uint} has been defined as {\tt unsigned int},
and {\tt LIMB} has been defined as {\tt unsigned short int}.

The {\tt MpiAdd} function adds long integer $B$ to long integer
$A$, and returns the carry. Both $A$ and $B$ are of length $n$,
that is, $A,B\in\{0..\rho^n-1\}$, and the carry can thus be {\tt
0} or {\tt 1}.

The {\tt MpiSub} function subtracts long integer $B$ from long
integer $A$, and returns the carry. Both $A$ and $B$ are of length
$n$, that is, $A,B\in\{0..\rho^n-1\}$, and the carry can thus be
{\tt 0} or {\tt -1}.

If $A=0$, then {\tt MpiNeg} does nothing but return {\tt 0} to
indicate $A$ is zero. If $A\neq 0$, then {\tt MpiNeg} replaces $A$
with $\rho^n-A$, and returns {\tt 1}.

The {\tt MpiAddC} function adds input carry $\kappa$ to long
integer $A$, and returns its output carry. This function is used
in carry processing.

\subsection{The Upper Level Group}

The upper level group of functions include the following
functions:

{\tt void MpiMul\_KR(uint $\rho$, LIMB* D, const LIMB* A, int n,
const LIMB* B, int m);}

{\tt void KRMpiMul(uint $\rho$, LIMB* D, const LIMB* A, int n,
const LIMB* B, int m);}

{\tt int KRMulTop(uint $\rho$, LIMB* D, const LIMB* A, const LIMB*
B, int n);}

{\tt int KRMul(uint $\rho$, LIMB* D, const LIMB* A0, const LIMB*
A1, const LIMB* B, int n);}

{\tt int KRMulB1(uint $\rho$, LIMB* D, const LIMB* A, const LIMB*
B, int n);}

{\tt int KRMulB2(uint $\rho$, LIMB* D, const LIMB* A0, const LIMB*
A1, const LIMB* B, int n);\/}

The {\tt MpiMul\_KR\/} is the highest level function. ``KR" stands
for ``Karatsuba and Roche". It computes the product of length-$n$
integer $A$ and length-$m$ integer $B$, and stores it in the $n+m$
limbs buffer starting at $D$. It calls {\tt KRMpiMul} to do the
actual work. The only thing that {\tt MpiMul\_KR\/} does itself is
to ensure that $n\geq m$ when calling {\tt KRMpiMul($\rho$, D, A,
n, B, m)}.

{\tt KRMpiMul} actually computes the product of length-$n$ integer
$A$ and length-$m$ integer $B$, and stores it in the $n+m$ limbs
buffer starting at $D$. It assumes $n\geq m$. When $m$ is smaller
than a threshold for using Karatsuba multiplication, it calls the
schoolbook multiplication algorithm. Otherwise, it computes the
product using the algorithm described in Section 2.4. It calls
{\tt KRMulTop} to do equal length additive multiplication.

{\tt KRMulTop} computes $AB+C\rho^n$, stores the result in the
$2n$ limbs buffer starting at $D$ and returns the carry, where the
long integer $C$ is stored beforehand in the first $n$ limbs of
$D$. It is a recursive function, and may call itself, {\tt KRMul},
, {\tt KRMulB1}, and {\tt KRMulB2}. {\tt KRMulTop} is based on a
simplified version of the algorithm described in Sections 2.2 and
2.3.

{\tt KRMul} computes $(A0-A1)B+C\rho^n$ and stores the result in
the $2n$ limbs buffer starting at $D$ and returns the carry, where
the long integer $C$ is stored beforehand in the first $n$ limbs
of $D$. It is a recursive function, and may call itself and {\tt
KRMulB2}. {\tt KRMul} is based on the full algorithm described in
Sections 2.2 and 2.3.

{\tt KRMulB1} serves the same purpose as {\tt KRMulTop}, using a
schoolbook algorithm. It is called when the length is less than
the Karatsuba threshold.

{\tt KRMulB2} serves the same purpose as {\tt KRMul}, using a
schoolbook algorithm. It is called when the length is less than
the Karatsuba threshold.

\section{Numerical Results}

The {\tt MpiMul\_KR} code has been thoroughly tested by the
present author. It was found to produce exactly the same results
as the schoolbook and the standard Karatsuba algorithms. To show
how efficient it is relative to the traditional algorithms, we
also made a numerical experiment. The experimental setting was as
follows:

\begin{itemize}
\item Average time per multiplication for lengths $100,
200,\cdots, 10000$ were measured for the three algorithms:
Schoolbook (SB), Karatsuba Standard (KS), space-efficient
Karatsuba or Karatsuba and Roche (KR). Only square multiplications
were considered in the experiment.

\item Before the experiment was done, we found that the Karatsuba
divide-and-conquer technique begins to gain benefit when $n\approx
100$, so we have made the Karatsuba threshold 128. This number was
chosen simply because it is a power of 2. There will be no
significant consequence if we choose a near but different
threshold.

\item The pairs of long integers for multiplication were generated
randomly. The number of pairs generated decreases as $n$
increases, since we wanted to finish the experiment in hours.

\item In coding Karatsuba standard, we incorporated some
speeding-up and memory saving techniques proposed in
\cite{bzmca,thome02}, so it is not purely conventional Karatsuba.
In our code each call of Karatsuba recursive function consumes
$2[n/2]$ limbs of heap memory.

\item Technical specifications of the computer running the
experiment were as follows: Windows 7 32-bit operating system,
Intel i3-2120 @3.30GHz CPU, and 4GB memory.
\end{itemize}

\begin{figure}[h]
\centering
\includegraphics[width=6in]{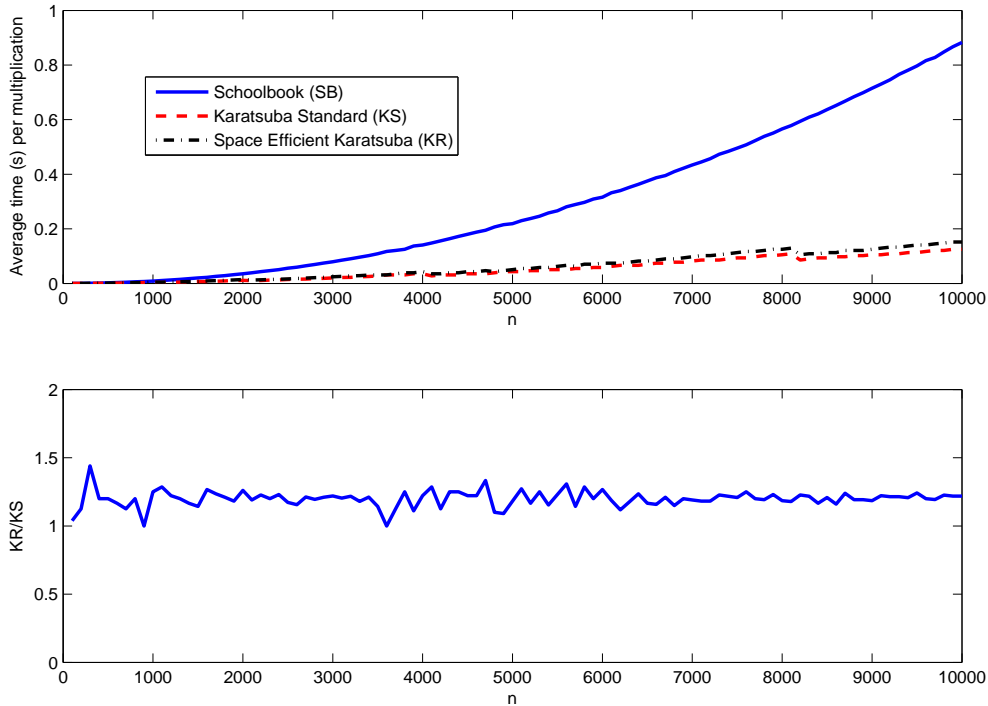}
\caption{Results of the numerical experiment comparing the
performances of SB, KS, and KR algorithms}
\end{figure}

The results are visualized in Figure 1. From the results data and
the figure we observe the following:

\begin{itemize}
\item The results seem to confirm the theoretical complexity
results. From $n=100$ to $n=10000$, for SB, the average execution
time grows 8678 times, close to the theoretical number 10000; for
KS, the average execution time grows 1382 times, and for KR,
average execution time grows 1622 times, both close to the
theoretical number $100^{\log_2 3}\approx 1479$.

\item The second chart shows KR takes roughly 20\% more time than
KS. Therefore, in terms of time efficiency, KR is roughly 20\%
less efficient than KS. This might be seen as a small price for
KR's avoiding of heap memory.
\end{itemize}

\section{Conclusion}

We have described in detail how to implement the space-efficient
Karatsuba algorithm for multiplication of long integers. Issues
resulting from carries, and other special issues for integers are
resolved. Numerical results show that the space-efficient
Karatsuba, while totally avoiding heap memory allocation, is
slightly less time efficient than the standard Karatsuba.

\bibliographystyle{unsrt}

\end{document}